\input amstex
\documentstyle{amsppt}
\magnification=\magstep1

\pageheight{9.0truein}
\pagewidth{6.5truein}

\TagsOnRight
\NoBlackBoxes

\hyphenation{uni-ser-ial}

\input xy
\xyoption{matrix}\xyoption{arrow}\xyoption{curve}

\long\def\ignore#1{#1}

\def\udotloopr#1{\ar@{{}.>} @'{@+{[0,0]+(-4,5)} @+{[0,0]+(0,10)}
  @+{[0,0]+(4,5)}} ^{#1}}

\def\seq{\mathrel{\widehat{=}}}
\def\la{{\Lambda}}
\def\lamod{\Lambda\text{-}\roman{mod}}
\def \len{\operatorname{length}} 
\def\detour{\mathrel{\wr\wr}}

\def\AA{{\Bbb A}}
\def\PP{{\Bbb P}}
\def\SS{{\Bbb S}}
\def\ZZ{{\Bbb Z}}
\def\NN{{\Bbb N}}
\def\Hom{\operatorname{Hom}}

\def\rank{\operatorname{rank}}
\def\stab{\operatorname{Stab}}
\def\rad{\operatorname{rad}}
\def\Og{{\Cal O}}
\def\GL{\operatorname{GL}}
\def\K{\operatorname{K}}

\def\Im{\operatorname{Im}}

\def\C{{\Cal C}}
\def\F{{\Cal F}}
\def\G{{\Cal G}}
\def\H{{\Cal H}}
\def\N{{\Cal N}}
\def\aut{\operatorname{Aut}}
\def\End{\operatorname{End}}

\def\T{{\Cal T}}
\def\uni{\operatorname{Uni}}

\def\moduni{\operatorname{Mod-Uni}}
\def\modla#1{{\operatorname{Mod}}_\Lambda^{#1}}
\def\dim{\operatorname{dim}}
\def\uniserdim{{\operatorname{uniserdim}}_{\SS}\Lambda}

\def\guni{{\Cal G}\operatorname{-Uni}}
\def\gunis{{\Cal G}\operatorname{-Uni}(\SS)}
\def\sub{\operatorname{Sub}^m_P}
\def\fac{\operatorname{Fac}^d_P}
\def\pr{\operatorname{pr}}
\def\BP{B^{m,d}_P}
\def\boldd{{\bold d}}
\def\algs{\operatorname{Alg}_s}
\def\algsmodr{{\operatorname{Alg}}_s{\operatorname{Mod}}^r}
\def\algssubm{{\operatorname{Alg}}_s{\operatorname{Sub}}^m}
\def\bmr{B^{m,r}}

\def\ARS{{\bf 1}}
\def\Bongiso{{\bf 2}}
\def\BongAdv{{\bf 3}}
\def\Bong{{\bf 4}}
\def\BongTrond{{\bf 5}}
\def\GeomII{{\bf 6}}
\def\Bor{{\bf 7}}
\def\Gab{{\bf 8}}
\def\GeomI{{\bf 9}}
\def\GeomIII{{\bf 10}}
\def\HuSk{{\bf 11}}
\def\Hum{{\bf 12}}
\def\Ki{{\bf 13}}
\def\Kraft{{\bf 14}}
\def\Kra{{\bf 15}}
\def\leb{{\bf 16}}
\def\Ros{{\bf 17}}

\topmatter
\title Varieties of uniserial representations IV. Kinship to geometric
quotients
\endtitle
\rightheadtext{Varieties of uniserial representations}

\author Klaus Bongartz and Birge Huisgen-Zimmermann\endauthor

\address FB Mathematik, Universit\"at Gesamthochschule Wuppertal, 42119
Wuppertal, Germany
\endaddress
\email bongartz\@math.uni-wuppertal.de\endemail

\address Department of Mathematics, University of California, Santa Barbara,
CA 93106\endaddress
\email birge\@math.ucsb.edu\endemail

 \thanks The research of the second author was partially supported by a
National Science Foundation grant.\endthanks


\abstract  Let $\la$ be a finite dimensional algebra over an algebraically closed field,
and $\SS$ a finite sequence of simple left $\la$-modules.  In
\cite{\GeomII, \GeomI}, quasiprojective algebraic varieties with accessible affine open
covers were introduced, for use in classifying the uniserial
representations of $\la$ having sequence $\SS$ of consecutive composition factors.  Our
principal objectives here are threefold:  One is to prove these varieties to be `good
approximations'  --  in a sense to be made precise  --  to geometric quotients of the
classical varieties $\moduni(\SS)$ parametrizing the pertinent uniserial representations,
modulo the usual conjugation action of the general linear group.  To some extent, this
fills the information gap left open by the frequent non-existence of such quotients.  A
second goal is that of facilitating the transfer of information among the `host'
varieties into which the considered uniserial varieties can be embedded.  These tools are
then applied towards the third objective, concerning the existence of geometric quotients:
We prove that $\moduni(\SS)$ has a geometric quotient by the $\GL$-action precisely when
the uniserial variety has a geometric quotient modulo a certain natural algebraic group
action, in which case the two quotients coincide.  Our main results are exploited in a
representation-theoretic context: Among other consequences, they yield a geometric
characterization of the algebras of finite uniserial type which supplements existing
descriptions, but is cleaner and more readily checkable.
\endabstract

\endtopmatter
\document

\head  Introduction \endhead

Our purpose is to study the uniserial representation theory of a finite dimensional algebra
$\la$ over an algebraically closed field $K$. This investigation falls
into natural parts, namely to describe the classes of uniserial representations sharing a
fixed sequence
$\SS$ of consecutive composition factors. One of the classical molds for the
classification problem is provided by the varieties $\moduni(\SS)$
consisting of those points in $\modla{d}$ which represent uniserial modules with
composition sequence $\SS$; here $d$ is the number of terms in 
$\SS$, and $\modla{d}$ is the traditional variety parametrizing all $d$-dimensional
$\la$-modules. In choosing this mold, one encounters major obstacles
on the road towards a classification of the uniserial modules  --  say, on the basis of
quiver and relations of $\la$  --  however: The open subvarieties $\moduni(\SS)$ of
$\modla{d}$ are very large, and even though geometric quotients by the natural
$\GL_d$-action (the orbits of which are in one-to-one correspondence with the isomorphism
types of representations) sometimes exist, such optimal situations are far from being
readily recognizable. In fact, the representation-theoretic information stored in these
varieties is encoded mainly in the
$\GL_d$-orbits, the geometry of which is notoriously difficult to access from a
presentation of
$\la$, in general; this is true even though the global varieties $\modla{d}$ are defined in
terms of such presentations. To meet these difficulties, finite collections
$V_\SS$ of smaller affine varieties parametrizing the considered uniserial representations
of $\la$ were introduced  in
\cite{\GeomI}. (For purposes of this introduction, we will refer to $V_\SS$ as though it
were a single variety.) These varieties provide a snug fit for the corresponding families
of uniserial representations  --  for more precision see below.   Moreover, they
have easily computable defining polynomials that stand in a transparent connection to
quiver and relations of $\la$, thus providing a viable bridge between quiver-presentations
and  uniserial representation theory.  In
\cite{\GeomII}, the  $V_{\SS}$ were shown to embed into certain Grassmannians
as locally closed subvarieties $\gunis$;  on the other hand, they are isomorphic to
certain  closed subvarieties of
the classical varieties $\moduni(\SS)$, albeit not to natural ones. 
 
Each of the above settings -- $\moduni(\SS)$, $V_\SS$, and $\gunis$ -- makes
certain aspects of the uniserial representation theory of
$\la$ comparatively easy to tackle;
first instances of this phenomenon can be found in \cite{\GeomII} and
\cite{\GeomIII}.  This motivates the wish for tools allowing for
smooth transfer of geometric information among the various `environments' into which the
considered uniserial varieties are embedded -- in particular, such shifts of perspective
are desirable as a great deal of general information on the varieties $\modla{d}$ is
available. Machinery for this purpose is provided in
Proposition C of Section 2, where we set up a general correspondence between
Grassmannian varieties of submodules of a projective module $P$ on one hand and 
 classical varieties of factor modules of
$P$ on the other.  

One of the points on our agenda is to
apply this philosophy of smoothing out technical obstacles through a change of scene 
to the following problem:  Namely, that of understanding the structure and size of the 
fibres of
the representation maps
$$\text{Considered parametrizing variety}  \ \longrightarrow \ \{\text{isom\. types of
uniserials with sequence\ } \SS\}.$$  
While this is a difficult task in the framework of the variety
$V_{\SS}$, as well as in that of its incarnation inside  
$\moduni(\SS) \subseteq \modla{d}$, it becomes easy  in
the Grassmannian setting $\gunis$. It turns out that the fibres of the representation
maps are always isomorphic to full affine spaces of small dimension; more precisely, if
$\SS = (S(0), \dots, S(l))$, where $S(0)$ is the top simple of the considered uniserial
modules, the fibre dimensions are bounded from above by the multiplicity of $S(0)$ in
$\SS$ diminished by $1$.  In addition, the fibres are closed subvarieties of $\gunis$,
which, by virtue of Proposition C, says that uniserial representations never have
proper uniserial degenerations (parts (2), (3) of Theorem A, Section 2).  On
learning about our results concerning the fibres, Le Bruyn partially re-proved them in
\cite{\leb}, with additional machinery, in the setting of the
appropriate Hesselink stratum of the nullcone of representations of a quiver modulo
relations.   In rough terms, the situation can be summarized by saying that
$\gunis$ is always a close approximation to a geometric quotient of $\moduni(\SS)$
modulo the $\GL_d$-action, in that the
$\GL_d$-orbits of $\moduni(\SS)$ are reduced to comparatively small closed
subvarieties of $\gunis$, all of which are geometrically harmless.  In other words, the
geometric  information encoded in the orbits of $\moduni(\SS)$ 
is shifted to the global geometry of $\gunis$.   Note however that, in
general, geometric quotients fail to exist.  In fact, the fibre dimension of the
representation map
$$\phi_{\SS}: \gunis \ \longrightarrow \ \{\text{isom\. types of uniserials with sequence\
} \SS\}$$ need not even be constant on the irreducible components (Example 1, Section 3). 
It is now natural to define the `uniserial dimension of
$\la$ at $\SS$' as the maximum of the differences `global dimension minus generic fibre
dimension', ranging over the irreducible components of $\gunis$.  As a consequence of the
close tie between quiver-presentations of
$\la$ and polynomials for $\gunis$  --  see part I of Section 3  --  the uniserial
dimension of
$\la$ at $\SS$ is readily available from such a presentation. 
This yields
manageable invariants of
$\la$ which play a crucial role in subsequent work on tame algebras
\cite{\HuSk}.  

Another of our aims is to decide when $\guni(\SS)$ actually {\it is\/} a geometric
quotient of $\moduni(\SS)$ by the $\GL_d$-action. The answer arises as a consequence of
the equivalence of two quotient problems: If $P$ is a projective
cover of $S(0)$, the variety $\guni(\SS)$ carries a natural $\aut_{\la}(P)$-translation
action and has a geometric quotient by this action if and only if $\moduni(\SS)$ has a
geometric quotient by the $\GL_d$-conjugation; in the positive case the two quotients are
isomorphic (part (1) of Theorem A).  The fact that
relating geometric quotients to other data on $\la$ is more manageable in the setting of
$\gunis$ than for the classical varieties $\moduni(\SS)$ is again attributable to the fact
that the former varieties are closer to such quotients to begin with.  Some special cases
of  our findings along this line overlap with work of King
\cite{\Ki}. We note that they also combine with
\cite{\GeomI, Theorem G} to show that each affine algebraic variety occurs as
a geometric quotient $\moduni(\SS) / \GL_d$ for suitable choices of $\la$ and
$\SS$, an indication of the immense complexity of the uniserial representation
theory of finite dimensional algebras in general.    

As already mentioned, one of the main virtues of the varieties $\gunis$ lies in
the fact that they have distinguished affine open covers which are tightly linked to quiver
and relations of $\la$.  Therefore they afford concrete insight into the uniserial
representation theory of the algebra $\la$, whenever the latter is given in this form. 
One of the focal points of investigations in this direction has been the problem of
characterizing the algebras which have only a finite number of uniserial modules up to
isomorphism, the algebras `of finite uniserial type' for short; it was posed by Auslander,
Reiten, and Smal\o\  in
\cite{\ARS, Problems (1) and (2) on p\. 411}.  Two previous articles addressing it
are by the first and second author, respectively.  In
\cite{\Bong}, an inductive characterization of finite uniserial type --  based on the Loewy
length of the algebra  --  is given, while
\cite{\GeomIII} provides combinatorial necessary and sufficient conditions in
terms of quiver and relations of $\la$, which, however, are separated by a slim gap
bridged by a system of linear equations. Two conjectures aiming at a
more manageable description were left open in \cite{\GeomIII}, namely: If
$\la$ has finite uniserial type, then (1)  all of the varieties
$\gunis$ are affine spaces, and (2) for each choice of $\SS$, there is at
most one uniserial module with this sequence of composition factors.  The
geometric results outlined above enable us to settle the first conjecture
in the positive and to thus supplement the confirmation of the second in
\cite{\Bong}.  Moreover, this second conjecture is re-obtained here with a geometric
argument.  The resulting information, in turn, leads to a new characterization of the
algebras of finite uniserial type which is cleaner and, from several viewpoints, more
satisfactory than existing ones (see Section 3, part II).  

Section 4, finally, contains an alternate proof for the fibre structures of the
representation maps, which yields the following interesting fact as a
by-product:  All endomorphism rings of uniserial $\la$-modules are commutative, a result
which, of course, relies heavily on algebraic closedness of the base field.

\head 1. Notation and tools \endhead

Throughout, $\Lambda$ will stand for a finite dimensional algebra over an
algebraically closed field $K$, and $J$ will be its Jacobson radical.  We
assume $\Lambda$ to be basic and may thus suppose that
$\Lambda = K\Gamma / I$ is a path algebra modulo relations, with underlying
quiver $\Gamma$.  It will be convenient to identify the vertices of $\Gamma$
with a full set of primitive idempotents of $\la$; all idempotents to which
we will refer will be taken from this set.  Moreover, we fix a sequence
$\SS = (S(0), \dots, S(l))$ of simples $S(i) =
\la e(i) / J e(i)$, where the $e(i)$ are primitive idempotents.  The uniserial modules $U$
of composition length
$l$ with
$J^{i}U/J^{i+1}U
\cong S(i)$ we call `uniserials with composition series $\SS$' for short, and to
the paths of length
$l$ passing through the corresponding sequence of vertices
$(e(0), \dots, e(l))$ in that order we will briefly refer as `paths through
$\SS$'. Finally, we call any element $x\in U\setminus JU$ with $e(0)x=x$ a {\it
top element} of $U$.   

The smoothest description of the varieties considered is in the framework of Grassmannians
as follows:	 Setting $m = \dim_K \la e(0) - (l+1)$, we denote the Grassmannian of all
$m$-dimensional $K$-subspaces of $\la e(0)$ by $\Cal G_m(\la e(0))$.  Then the
quasiprojective subvariety $\gunis \subseteq \G_m(\la e(0))$ consists of those
$m$-dimensional sub{\it modules} $C$ of $\la e(0)$ for which
$\la e(0)/C$ is uniserial with composition series $\SS$. As we saw in
\cite{\GeomII},
$\gunis$ is a locally closed subset of $\G_m(\la e(0))$ having a distinguished affine
open cover as follows:  Given a uniserial module $U$ with composition series $\SS$, each
path through $\SS$ which does not annihilate $U$ is  called a {\it
mast} of $U$.  Letting $p$ be a path through $\SS$ and $\guni(p)$ the set of points $C \in
\gunis$ such that
$\la e(0)/ C$ has mast $p$, we obtain a family   
$$(\guni(p))_{p\ \text{a path through}\ \SS}$$
of open affine subsets exhausting $\gunis$.   Clearly, $\gunis$ comes 
equipped with a natural surjection $\phi_{\SS}$ onto the set of isomorphism classes of
uniserial modules with composition series
$\SS$ which maps a point $C \in \gunis$ to the class of $\la e(0) / C$.  Moreover, it is
obvious that the subsets 
$\guni(p) \subseteq \gunis$ are unions of fibres
of $\phi_{\SS}$.  Therefore, we can safely
narrow our focus to the affine situation in exploring the fibres. We
denote the restriction of the map $\phi_\SS$ to $\guni(p)$ by $\phi_p$.

In \cite{\GeomII} it was also shown that the affine varieties $\guni(p)$ have
several isomorphic incarnations within alternate settings which frequently offer
advantages in solving specific problems.  Since, in the sequel, we will move freely
back and forth among these varieties, depending on which viewpoint offers an edge, we
will briefly recall their definition.

The first of the alternate incarnations, introduced in \cite{\GeomI}, is combinatorial in
nature and requires a bit of terminology for convenient communication of the
relevant data.  Fix a path $p$ through
$\SS$, and denote by $p_0, \dots, p_l$ its {\it right subpaths}; this means that $p =
q_i p_i$ for suitable paths $q_i$, our convention being to write $q_i p_i$ for
`$q_i$ after $p_i$'.  We assume the paths $p_i$ to be ordered in terms of
increasing lengths, i.e., $\len(p_i) = i$.  A path
$v \in K\Gamma$ starting in the vertex $e(0)$ is called a {\it route} on $p$ if
the sequence of vertices through which $v$ successively passes is a subsequence
of $(e(0), \dots, e(l))$.  Examples of routes are paths $\alpha p_m$ based
on `detours' $(\alpha, p_m)$, as follows: We say that a pair $(\alpha,p_m)$,
combining an arrow
$\alpha$ with a right subpath $p_m$ of $p$, is a {\it detour} on
$p$ and write $(\alpha,p_m) \detour p$, in case $\alpha p_m$ fails to be a
right subpath of $p$, while there exists a right subpath $p_s$ longer than $p_m$
which has the same terminal point as $\alpha$. The set of all such indices $s$ is denoted
by
$I(\alpha, p_m)$.  Next we observe that, given any uniserial module $U$ with mast $p$ and
top element $x$, we have $\alpha p_m x = \sum_{i \in I(\alpha,p_m)} k_i(\alpha,p_m) p_i x$
for unique scalars $k_i(\alpha, p_m)$.  The prominent role played by detours is
due to the obvious fact that it suffices to record the effect of
multiplying top elements of $U$ by detours, in order to pin down the isomorphism
type of
$U$.  Reflecting this, the affine variety $V_p$ lives in
$\Bbb A^N$, where $N = \sum_{(\alpha,u) \detour p} |I(\alpha,u)|$; it
consists of the points $k$ which show up as coordinate
strings $(k_i(\alpha,u))_{i \in I(\alpha,u), (\alpha, u) \detour p}$ of
uniserials with mast $p$ relative to arbitrary choices of top elements.  (For the fact that
this set of points in $\AA^N$ actually is an affine algebraic variety, see \cite{\GeomI}.) 
An isomorphism $\psi_p: V_p
\rightarrow \guni(p)$ is afforded by the assignment 
$$k = (k_i(\alpha,u)) \mapsto  \biggl( \sum_{(\alpha,u) \detour p}
\la \bigl(\alpha u \ -
\sum_{i \in I(\alpha,u)} k_i(\alpha,u)p_i \bigr) \ \  + \sum_{q \ \text{not a
route on}\ p} \la qe(0) \biggr)$$
(see \cite{\GeomII}).  If $\Phi_p$ denotes the surjection from $V_p$ onto
the set of isomorphism classes of uniserials with mast $p$, which takes $k$ to $\la e(0)
/\psi_p(k)$, we clearly obtain the following commutative triangle:  

\ignore{
$$\xymatrixcolsep{4pc}\xymatrixrowsep{1pc}
\xymatrix{
 V_p \ar[dr]^-{\Phi_p} \ar[dd]_{\psi_p}^{\cong}\\
 & \save+<23ex,0ex> \drop{\txt{ \{iso types of uniserials in $\lamod$ with mast
$p$\} }} \restore\\
\guni(p) \ar[ur]_-{\phi_p} }$$ }

\noindent

There is another helpful way of looking at the affine varieties $\guni(p) \cong
V_p$.  Occasionally, it is convenient to identify the affine variety
$V_p$ (alias $\guni(p)$) with a subvariety, labeled $\uni(p)$, of the classical variety
$\modla{l+1}$ of bounden $(l+1)$-dimensional representations of $\la$; this
identification is justified in \cite{\GeomII, Section 3}.  Recall that, given a natural
number $d$, the variety 
$\modla{d}$ of bounden $d$-dimensional representations of
$\Lambda = K\Gamma / I$ is defined as follows:  If $\Gamma^* = \{ \text{vertices of}\
\Gamma\} \cup \{
\text{arrows in}\ \Gamma\}$, then
$\modla{d}$ consists of those points $x = (\alpha(x)) \in \prod_{\alpha \in
\Gamma^*} M_d(K)$, the components of
which satisfy the relations in $I$.  More precisely, $\uni(p)$ lives inside
the open subvariety $\moduni(p)$ of $\modla{l+1}$ containing precisely those
points which correspond to uniserial modules with mast $p$. 
Namely, if $p = \alpha_l \alpha_{l-1}
\cdots \alpha_1$, where the $\alpha_i$ are arrows, then $\uni(p)$ consists of all
those points $$x = (\alpha(x))_{\alpha\ \in \ \Gamma^*} \in \moduni(p),$$
for which the $i$-th column of the $(l+1) \times (l+1)$ matrix
$\alpha_i(x)$ equals the $(i+1)$-st canonical basis vector for $K^{l+1}$, $1 \le
i \le l$;  note that the modules corresponding to these points are
automatically uniserial with mast $p$.  The canonical isomorphism
$V_p \rightarrow \uni(p)$ actually identifies the map $\Phi_p$ with the
restriction of the canonical map $R$ from $\modla{l+1}$ to the isomorphism types of
left
$\la$-modules; so in particular, fibres of $\Phi_p$ are carried to fibres of $R$ by the
mentioned isomorphism of varieties.  Note that
$\uni(p)$ is no longer stable under the conjugation action of $\GL_{l+1}$; the
closure of $\uni(p)$ under $\GL_{l+1}$-conjugation coincides with the full
subvariety $\moduni(p)$ of $\modla{l+1}$.  In fact, the natural group action on
$\guni(p)$ $(\cong V_p \cong \uni(p))$, which will provide our main operative tool,
does not translate into any conjugation action of a matrix group on
$\uni(p)$ in general (see Example 1).

We conclude this sketch of background information with an overview of the
mentioned perspectives, including brief
descriptions of the relevant varieties. Suppose that $p$ is a path in $\Gamma$ of length
$l$.
\smallskip

\noindent $\bullet$ $\guni(p)$ is the subvariety of the Grassmannian $\G_{\dim_K \la e(0)
-(l+1)}(\la e(0))$ consisting of all
$\la$-submodules $C\subseteq \la e(0)$ such that $\la e(0)/C$ is uniserial with mast $p$;
moreover,
$\gunis$ is the union of the $\guni(p)$, where $p$ traces all paths through $\SS$.
\smallskip

\noindent $\bullet$ $V_p$ is the set of all coordinate strings glued together from all the
coordinate vectors of detours inside uniserials with mast $p$.
\smallskip

\noindent $\bullet$ $\moduni(p)$ is the open subvariety of $\modla{l+1}$ consisting of the
points going with uniserials modules that have mast $p$, and $\moduni(\SS)$ is their union
as $p$ again traces the paths through $\SS$.
\smallskip

\noindent $\bullet$ The $\uni(p)$ are the subvarieties of the corresponding $\moduni(p)$
introduced in the preceding paragraph, and $\uni(\SS)$ is their union.
\smallskip

These varieties are related as shown in the following diagram:

\ignore{
$$\xymatrixcolsep{4pc}\xymatrixrowsep{3pc}

\xymatrix{
\guni(p) \ar[r]^-{\cong} \ar[d]_-{\subseteq} &V_p \ar[r]^-{\cong} &\uni(p)
\ar[r]^-{\subseteq} \ar[d]_-{\subseteq} &\moduni(p)
\ar[r]^-{\subseteq} \ar[d]_-{\subseteq} &\modla{l+1} \ar@{=}[d]\\
\gunis  &&\uni(\SS) \ar[r]^-{\subseteq} &\moduni(\SS) \ar[r]^-{\subseteq}
&\modla{l+1} 
}$$
}

\noindent Note that the
first row of our diagram consists entirely of affine varieties.  As for the
second row:  On one hand, 
$\gunis$ often fails to be affine, while on the other hand, $\uni(\SS)$ is
a closed subset of the affine variety $\modla{l+1}$.  In particular, this shows
that the blank in the second row cannot be filled by an isomorphism in general.

The legitimacy of moving from one incarnation of $\guni(p)$ to another, as
convenience dictates, is guaranteed by the following commutative diagram which
ties all of the representation maps together:

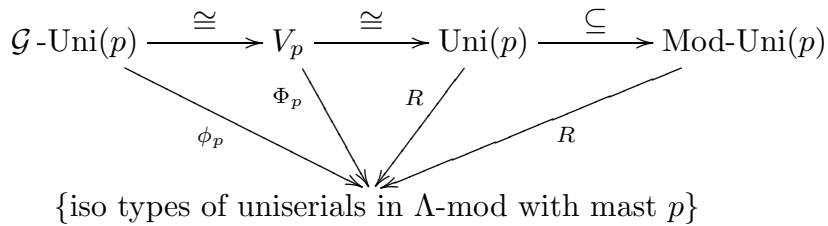
\ignore{
$$\xymatrixcolsep{1.25pc}\xymatrixrowsep{3pc}
\xymatrix{ 
\guni(p) \ar[rr]^-{\displaystyle\cong} \ar[drrr]_(0.5){\phi_p} &&V_p
\ar[rr]^-{\displaystyle\cong} \ar[dr]_(0.3){\Phi_p} &&\uni(p)
\ar[rr]^-{\displaystyle\subseteq} \ar[dl]_(0.45){R} &&\moduni(p)
\ar[dlll]^(0.5){R}\\
 &&& \save+<0ex,-1.2ex> \drop{\txt{ \{iso types of
uniserials in $\lamod$ with mast $p$\} }} \restore
}$$
}

\noindent Here $R$ denotes the canonical map from
$\modla{l+1}$ to the set of isomorphism types of $(l+1)$-dimensional left
$\la$-modules.

\head 2. Fibre structure and transfer of information
\endhead

Let $G$ be a linear algebraic group acting morphically on a variety $X$. We
will use the strongest notion of a quotient of
$X$ by $G$, which is as follows:  A surjective morphism $\phi: X
\rightarrow Y$ of varieties is a {\it geometric quotient} of $X$  modulo $G$
if $\phi$ is  open, the fibres of $\phi$ coincide with the orbits of $G$,
and, for any open subset $U
\subseteq Y$,  the comorphism $\phi^0$ of $\phi$ induces an isomorphism from
the ring $\Og_Y(U)$ of regular functions on $U$ to the subring of 
$\Og_X(\phi^{-1}(U))$  consisting of those regular functions which are constant on
the $G$-orbits of $\phi^{-1}(U)$.  Any such geometric quotient is a {\it
categorical} one, in the sense that each morphism $X \rightarrow Z$ which is
constant on the $G$-orbits of $X$ factors uniquely through $\phi$;  in
particular, geometric quotients are unique up to isomorphism in case of
existence.  

Let us focus on a sequence $\SS$ of simple left $\la$-modules as before and
abbreviate the idempotent $e(0)$ corresponding to the top simple by
$e$.  We start by introducing a left action of the algebraic group
$\aut_{\la}(\la e)$ on $\gunis$ as follows: 
Given a point $C \in \gunis$, we define $gC$ to be $g(C)$.  Since clearly
$g$ induces an isomorphism of $\la$-modules $\la e / C
\cong \la e / gC$, this yields a well-defined morphic action which leaves
the fibres of $\phi_{\SS}$ invariant.  The action is transitive on
the fibres of $\phi_{\SS}$:  Indeed, due to the fact that $\la e$ is
projective and local, any isomorphism
$\la e / C \cong \la e / C'$ is induced by an automorphism of $\la e$ taking
$C$ to $C'$.  Moreover, since the isomorphic uniserials $\la e/ C$ and
$\la e/ g(C)$ have the same masts, the affine patches
$\guni(p)$ are unions of orbits.  In other words: For
each path $p$ through $\SS$, the $\aut_{\la}(\la e)$-action on $\gunis$
restricts to an action on $\guni(p)$, the orbits of which coincide with the
fibres of the canonical map $\phi_p$.   

The unipotent group $G$ which will work for us is the unipotent radical
of $\aut_{\la}(\la e)$.  A typical element of $G$ thus corresponds to right
multiplication of $\la e$ by a local unit from the set
$\{g = e + g' \mid  g' \in eJe\}$.   Clearly, $\aut_{\la}(\la e) \cong K^*
\times G$, and since $K^*$ stabilizes all elements of
$\gunis$, the orbits of the $G$-action are still identical with the
fibres of the canonical maps.  In particular, $G$ stabilizes all of the affine
patches $\guni(p)$.

This essentially yields the second part of our main theorem.

\proclaim{Theorem A} Let $\SS= (S(0),\dots,S(l))$ be a sequence of simple left
$\la$-modules, and let $\mu(\SS)$ denote the multiplicity of $S(0)$ in $\SS$.
\smallskip

\noindent{\rm(1)}  The classical variety $\moduni(\SS)$ has a
geometric quotient by the standard $\GL_d$-con\-ju\-ga\-tion precisely when $\gunis$ has
a geometric quotient by its $\aut_{\la}(\la e)$-action.  In the positive
case, the two geometric quotients are isomorphic.  Moreover, a geometric quotient of
$\moduni(\SS)$ by $\GL_d$ exists and coincides with $\gunis$ if and only if $\mu(\SS)$
 equals the $K$-dim\-en\-sion of the endomorphism ring of any
uniserial module with composition series $\SS$. 
\smallskip
 
\noindent {\rm(2)}  There is a morphic action of the unipotent group $G$ on
$\gunis$, the orbits of which coincide with the fibres of $\phi_{\SS}$.  In
particular, the fibres are closed subvarieties, and the fibre
dimension is generically constant on the irreducible components of $\gunis$.
\smallskip

\noindent {\rm(3)} Each fibre of $\phi_{\SS}$ is a
homogeneous $G$-space under this action.  In particular, each fibre is isomorphic
to a full affine space.  In fact, if $C \in \gunis$, then
$G.C = \phi_{\SS}^{-1}\phi_{\SS}(C) \cong \AA^{m(C)}$, where
$$m(C) =  \mu(\SS) - \dim_K\End_\la(\la e / C).$$
\endproclaim

Of course, analogous statements concerning the fibres hold for the canonical surjections
$\Phi_p$ from the affine varieties $V_p$ onto the set of isomorphism types of uniserial
modules with mast $p$.  It is this rendering which makes the
quantities arising in the theorem concretely accessible, an aspect addressed
in part I of Section 3. 

We will briefly discuss two immediate consequences of Theorem A and then follow with a
proof, which will occupy the remainder of this section.

Parts (2) and (3) of Theorem A prompt us to define the uniserial dimension of
$\la$ at the sequence $\SS$ of simple left $\la$-modules as follows.
 
\definition{Definition}  The {\it
uniserial dimension of $\la$ at $\SS$}, denoted $\uniserdim$, will be $-1$ in
case
$\gunis = \varnothing$; otherwise,
$\uniserdim$ will stand for the supremum of the following
differences:
$\dim \C$ minus the generic fibre dimension of $\phi_{\SS}$ on $\C$, where $\C$
runs through the irreducible components of $\gunis$.
\enddefinition

The uniserial dimension of $\la$ at $\SS$ is an
isomorphism invariant of $\la$.  Note that it, too, is readily computed by
way of the $V_p$: Since the subvarieties $\guni(p) \cong V_p$ form an open
cover of $\gunis$ the irreducible components of $\gunis$ are birationally
equivalent to those of the $V_p$'s, where $p$ traces the paths through $\SS$.

Let us consider the extreme case $\uniserdim = 0$.  As is to be expected, it
occurs precisely when there is a finite positive number of uniserial left
$\la$-modules with composition series $\SS$:  Indeed, that finiteness of this
number implies
$\uniserdim \le 0$ follows directly from the definition.  Conversely, if $\uniserdim = 0$,
the number of fibres of $\phi_{\SS}$ is indeed finite, for closedness of the fibres in
$\gunis$ guarantees that each irreducible component consists of a single
fibre.

As for part (1) of Theorem A:
Clearly, the necessary and sufficient condition for $\gunis$ to be a geometric quotient of
$\moduni(\SS)$ is satisfied if
$\mu(\SS) =1$, or if
$\mu(\SS) = 2$ and $S(0) \cong S(l)$.  In general, a presentation of
$\la$ in terms of quiver and relations permits us to check whether the condition
is satisfied (cf the computational remarks in Section 3). 

It was pointed out to us by W. Crawley-Boevey that the
existence of a geometric quotient of $\uni(\SS)$ modulo the $\GL_d$-action
in case $\mu(\SS) =1$ can also be derived from King's work on moduli
spaces \cite{\Ki} as follows:  If the quiver $\Gamma$ has $n$ vertices
$e_1, \dots, e_n$ with $e_1 = e(0)$, consider the homomorphism $\Theta:
\K_0(\lamod) \rightarrow \ZZ$ sending $(s_1, \dots, s_n)$ to $- ls_1 +
\sum_{1 \le i \le n} s_i$, where $l$ is the length of the sequence $\SS =
(S(0), \dots, S(l))$.  If $\mu(\SS) = 1$, then clearly all uniserial
modules with series $\SS$ are $\Theta$-stable and \cite{\Ki, Proposition 5.3}
applies.
\medskip 

Combining part (1) of Theorem A with
\cite{\GeomI, Theorem G}, finally, we obtain the following

\proclaim{Corollary B}  Every affine algebraic variety $V$ over $K$ occurs as a
geometric quotient of $\moduni(\SS)$ by the $\GL_d$-action for some
finite dimensional algebra and some sequence $\SS$ of simple modules.  One can
even arrange for isomorphisms 
$V \cong \moduni(\SS) / \GL_d  \cong \gunis$. \qed
\endproclaim

\medskip
In rough terms, the proof of part 1 of Theorem A (below) will depart from the following
obvious representation-theoretic correspondence:  Namely, given a projective
module
$P$ of $K$-dimension $r$, and natural numbers $m$ and $d$ with $r = m + d$, there is
a bijection between the $\aut_{\la}(P)$-orbits of the $m$-dimensional submodules
of $P$ and the isomorphism types of $d$-dimensional factor modules of $P$; it is
induced by $U \rightarrow P/U$.  We will translate this bijection into
a correspondence between sets contained in a subvariety of a certain Grassmannian
on one hand, and sets contained in a subvariety of $\modla{d}$ on the other
(Proposition C below).  The former subvariety is the closed subset $\sub$ of
$\G_m(K^r)$ consisting of all
$m$-dimensional submodules of $P$; here $\G_m(K^r)$ denotes the Grassmannian of
$m$-dimensional subspaces of $K^r$.  The subvariety $\fac$ of $\modla{d}$ which
we will consider consists of those points which correspond to
$d$-dimensional homomorphic images of $P$.  It is an open subvariety of
$\modla{d}$:  Indeed, a $d$-dimensional module $M$ belongs to $\fac$ if and only if
the multiplicity of any simple module $S$ in $M/ JM$ is smaller than or equal
to the multiplicity $c(S)$ of $S$ in $P / \rad P$.  But each of these
multiplicities is given by the integer-valued, upper semi-continuous function
$X \mapsto \dim\Hom_{\la}(X,S)$, and so the set
$\fac$, arising as a finite intersection of preimages $(- \infty, c(S) +1)$,
is indeed open.  We observe that $\aut_\la(P)$ naturally acts on $\sub$, while we
have the standard conjugation action of $\GL_d$ on $\fac$.  

In order to relate the projective variety $\sub$ to the quasi-affine variety
$\fac$  --  or, more precisely, to relate the $\aut_\la(P)$-stable subsets of
$\sub$ to the $\GL_d$-stable subsets of $\fac$  --  we introduce an intermediary,
namely the following closed subvariety $\BP$ of
$\GL_r$:  It contains those matrices in $\GL_r$ whose first
$m$ columns generate a $\la$-submodule of $P$.  In other words, if 
$p =(\alpha(p))_{\alpha \in \Gamma^*}$ is a point in
$\modla{n}$ representing $P$, the spaces
generated by the first $m$ columns of the matrices in $\BP$ are
invariant under left multiplication by all the matrices
$\alpha(p)$.  Clearly, $\aut_\la(P)$ acts on $\BP$.  So does the 
subgroup $\H$ of $\GL_r$, consisting of the upper triangular block
matrices of the form $\left[\smallmatrix h_1 & h_2 \\ 0 &
h_3\endsmallmatrix \right]$, where $h_1$ and
$h_3$ are invertible $m \times m$ and $d \times d$-matrices, respectively. 
  In summary, we obtain a right
action of the group $\aut_\la(P) \times \H$ on $\BP$, given by $b(f,h) =
f^{-1}bh$.  Clearly, suitable restriction of the right $\H$-action on $\BP$
provides us with a $\GL_d$-action on $\BP$ as well.  Whenever called for, we
identify $\GL_d$ with its canonical copy inside $\H$. 

The next step is to set up a diagram of morphisms 
$$\sub @<\rho<< \BP @>\sigma >> \fac,$$ where $\rho$ and $\sigma$ are
sufficiently well-behaved to afford the transit of geometric information, as we
shift subsets from $\sub$ to $\fac$ by applying
$\sigma$ to preimages under $\rho$, and vice versa.  In particular, we
want $\rho$ and $\sigma$ to be equivariant relative to the actions of
$\aut_\la(P)$ and $\GL_d$, respectively.  Moreover, we will construct both maps
so as to have local sections, which will at least give us local morphisms between
the two varieties flanking $\BP$. 

Here are the details of the setup: If we let $\rho: \BP \rightarrow \sub$ be the
map sending a matrix $b \in \BP$ to the space generated by its first $m$ columns,
then $\rho$ is clearly an $\aut_{\la}(P)$-equivariant morphism.  To define
$\sigma : \BP \rightarrow \fac$, we start by introducing the map $\pr: M_r(K)
\rightarrow M_d(K)$ which assigns to any $r \times r$-matrix its lower right $d
\times d$-block.  Moreover, we recall that the coordinates $\alpha(p)$ of the
point $p$ representing our projective module $P$ are labeled by the elements of
$\alpha \in \Gamma^* = \{e_1, \dots, e_n\} \cup \{$arrows in $\Gamma\}$.  Noting
that, for
$b \in \BP$ and any $\alpha \in \Gamma^*$, the matrix $b^{-1}\alpha(p) b$
is a block matrix of the form
$$\left( \smallmatrix * & * \\ 0 & \alpha(v) \endsmallmatrix \right)$$
with $\alpha(v) = \pr(b^{-1}\alpha(p) b)$, we obtain a point $v =
(\alpha(v)) \allowmathbreak \in \modla{d}$ which represents $P/ \rho(b)$  (by a
slight abuse of notation, we identify the point $\rho(b)$ of $\sub$ with the
corresponding $\la$-module).  Now we define
$\sigma$ to be the morphism 
$$b \mapsto  (\pr (b^{-1}\alpha(p) b))_{\alpha\in \Gamma^*}.$$ 
It is easily seen that $\sigma$ is, in fact, a surjective morphism which is
equivariant relative to the described operations of
$\GL_d$ on domain and codomain.  Equivariance is just a slice of the good
behavior of $\rho$ and $\sigma$, however.  

\proclaim{Lemma 1} The morphism $\rho : \BP \longrightarrow \sub$ is an
$\aut_{\la}(P)$-equivariant principal right $\H$-bundle. It is the
geometric quotient of $\BP$ by the action of $\H$. \endproclaim

\demo{Proof} Clearly, $\rho$ is the restriction of the $\GL_r$-invariant
morphism $\tau: \GL_r \rightarrow \G_m(K^r)$, again assigning to an
invertible matrix the space generated by its first $m$ columns.  Since the
canonical right action of $\H$ on $\GL_r$ restricts to the $\H$-action on
$\BP$, it suffices to prove that $\tau$ is a
$\GL_r$-equivariant principal right
$\H$-bundle. 
 
Equivariance being clear, we wish to provide a suitable open covering for
$\G_m(K^r)$.  To that end, let
$V$ be any subspace of dimension $d$ of $K^r$, and $\G_V$ the set of
points in
$\G_m(K^r)$ which complement $V$ in $K^r$.  The set $\G_V$ is an open affine
subset of $\G_m(K^r)$:  Indeed, choose a basis $b_{m+1}, \dots, b_r$ of
$V$, supplement it to a basis $b_1, \dots, b_r$ of $K^r$, and let $\Cal B$
be the corresponding basis of $\bigwedge^m(K^r)$.  Then $\G_V$ consists
precisely of those points  which have nonzero $(b_1 \wedge \dots\wedge
b_m)$-coefficient with respect to $\Cal B$.

To find suitable isomorphisms $\tau^{-1}(\G_V) \rightarrow \G_V \times
\H$, we observe that, due to the $\GL_r$-equivariance of $\tau$, we may
restrict our attention to the subspace $V\subseteq K^r$ which is generated
by the last $d$ canonical basis vectors of $K^r$.  Then each element
$x \in \tau^{-1}(\G_V)$ can be uniquely factored in the form $c\cdot h$ with
$h \in \H$ and $c$ a lower triangular matrix of the form $\left[\smallmatrix
c_1 & 0 \\ c_2 & c_3\endsmallmatrix \right]$, where $c_1, c_3$ are the
identity matrices of sizes
$m \times m$ and $d \times d$, respectively. The assignment $x \mapsto
(\tau(x),h)$ yields an isomorphism $\tau^{-1}(\G_V) \rightarrow \G_V \times
\H$ of varieties, and as $V$ varies, these
isomorphisms satisfy the required compatibility conditions.  

In particular, the morphism $\rho$ has local sections, and hence the final
statement now follows from
\cite{\BongTrond, Lemma 5.5}.
\qed \enddemo

There is an alternate guise of the
$\GL_d$-conjugation on $\fac$ which will be useful in the proof of Lemma
2: Namely, if we define a right action of
$\aut_{\la}(P) \times \H$ on $\fac$ by
$v (f,h)  =  (\pr h)^{-1} v (\pr h)$, then the morphism
$\sigma$ becomes $\aut_{\la}(P) \times \H$-equivariant, since the fibres of
$\sigma$ are stable under $\aut_\la(P)$.  Noting that the map $\pr :
M_r(K) \rightarrow M_d(K)$  restricts to a group homomorphism
$\pr|_{\H} : \H\rightarrow \GL_d$, and denoting its kernel by $N$, we obtain

\proclaim{Lemma 2} The morphism
$\sigma: \BP \rightarrow \fac$ is $\GL_d$-equivariant and smooth.  It is the
geometric quotient of $\BP$ by the $\aut_{\la}(P) \times N$-action.
\endproclaim

\demo{Proof}  For our analysis of $\sigma$ it will be helpful to
`expand' the picture and look at the following isomorphic copy $Z$ of the
variety $\BP$:  It is located as a locally closed subvariety inside the product 
$$\fac \times K^{d \times r} \times \BP,$$
namely as the image under the morphism $\zeta$ sending any point $b \in \BP$
to
$(\sigma(b), [0\ E]b^{-1}, b)$, where $E$ is the $d \times d$ identity
matrix.  One checks that the map
$\psi =  [0\ E]b^{-1}$ satisfies $\psi \alpha(p) = \alpha(\sigma(b))$ for
$\alpha \in \Gamma^*$ and thus concludes that $\psi$ is a $\la$-epimorphism
from $P$ to $P/\rho(b)$.  Clearly, $\zeta: \BP
\rightarrow Z =
\Im(\zeta)$ is a bijective morphism, the inverse of which is just the projection
onto the third component, whence $\zeta$ is an isomorphism.  The map $\zeta$ is
even an
$\aut_{\la}(P) \times \H$-equivariant isomorphism if one equips $Z$ with the
right $\aut_{\la}(P) \times \H$-action given by 
$$(\sigma(b), [0\ E]b^{-1}, b) \cdot (h,f) = \bigl((\pr h)^{-1} \sigma(b) (\pr
h), [0\ E](bh)^{-1} f, f^{-1}bh \bigr)$$ 
for $b \in \BP$,  $f \in \aut_{\la}(P)$ and $h \in \H$.  In verifying our
assertions, we are hence free to replace
$\sigma$ by the projection $\pi: Z \rightarrow \fac$ onto the first component.
 
First we check that the fibres of $\pi$ coincide with the $\aut_{\la}(P)
\times N$-orbits of $Z$ under the specified action.  The group $N$ being equal to
$$\left\{ \left[\smallmatrix h_1 & h_2 \\ 0 & E\endsmallmatrix \right]
\mid h_1
\in \GL_m, h_2 \in K^{m \times d} \right\},$$
it is clear that $\sigma$, and hence also $\pi$, is constant on these orbits. 
Conversely, consider two points in $Z$ with the same image under $\pi$, say $(v,
\psi, b)$ and $(v, \psi', b')$.  Then $\psi$ and $\psi'$ both
represent epimorphisms $P \rightarrow P/\rho(b)$, where the factor module
$P/\rho(b)$ represents the point $v \in \fac$, and by
factoring each of the two maps into a projective cover of $P/U$ and a trivial
component, we see that $\psi' = \psi f$ for a suitable $\la$-automorphism
$f$ of $P$; in other words, $\psi'$ belongs to the
$\aut_{\la}(P)$-orbit of $\psi$, and it is harmless to assume $\psi =
\psi'$.  But this implies that $(b')^{-1} = h b^{-1}$ for a
suitable element $h \in N$ and thus forces $b'$ into the $N$-orbit
of $b$.

For the first assertion of the lemma, it now suffices to show that
$\pi$ has local sections (see \cite{\BongTrond, Lemma 5.5}).  We prove this in
tandem with smoothness, by factoring $\pi$ in the form $\pi = \pi''\pi'$,
the first map being the projection
$$\ pi' : Z \rightarrow X, \ \ \ (v,\psi,b) \mapsto (v,\psi),$$  and the second 
$$\pi'' : X \rightarrow \fac, \ \ \ (v,\psi) \mapsto v;$$
here $X$ is the variety consisting of all pairs $(v,\psi) \in \fac \times K^{d
\times r}$ such that $\rank (\psi) = d$ and $\psi \alpha(p) = \alpha(v) \psi$
for $\alpha \in \Gamma^*$.  In other words, $\psi$ is a
$\la$-epimorphism $P \rightarrow P/U$, where $P/U \in \lamod$ corresponds
to $v$. 

To see that $\pi'$ has local sections, let $(v,\psi)$ be a point in $X$, and
write $\psi = (\psi_1, \psi_2)$, where $\psi_1 \in K^{d \times m}$ and $\psi_2
\in K^{d \times d} = M_d(K)$.  Since $\rank(\psi) = d$, we can find $g \in
\GL_r$ such that $(\psi g)_2 \in \GL_d$.  Pick an open
neighborhood $\N$ of $(v,\psi)$ in $X$, with the property that $(\psi' g)_2 \in
\GL_d$ for all points $(v',\psi') \in \N$, and check that the morphism $\N
\rightarrow Z$, given by $(v',\psi') \mapsto (v',\psi', s(v',\psi')),$ with
$$s(v',\psi') \ \ \ = \ \ \ g\left[\smallmatrix E_m & 0 \\
-(\psi'g)_2^{-1}(\psi'g)_1 & (\psi'g)_2 \endsmallmatrix \right],$$
is a section of $\pi'$ over $\N$.  Since all of the maps $\psi'$ are
surjective, Lemma 3.4 of \cite{\BongTrond}, and an obvious `affine shift'
thereof, yields the following two bundles $Z_0$ and $Z_1$ with base $\N$: 
Namely, the vector bundle $Z_1$ consisting of all triples $(v',\psi',c)$ with
$(v',\psi') \in \N$ and $c \in M_r(K)$ satisfying $\psi'c = 0$, and the
affine bundle $Z_0$ consisting of the triples $(v',\psi',b)$ with
$(v',\psi') \in \N$ and $b \in M_r(K)$ satisfying $\psi'b = [0 E]$.  We note that
$Z_0$ contains $Z$ as an open subset:  Indeed, given $(v',\psi',b) \in Z_0$, we
have
$b \in \BP$ if and only if $b \in \GL_r$, since the equality $\psi'b = [0 E]$
forces the first $m$ columns of $b$ to span the kernel of the $\la$-homomorphism
$\psi'$.  Moreover, we observe that $Z_1$ and $Z_0$ are isomorphic bundles via
the assignment $(v',\psi',c) \mapsto (v',\psi', c + s(v',\psi'))$.  Hence $\pi'$
is smooth, being the restriction of the vector bundle projection $Y \rightarrow
\N$ to the open subvariety $Z$ of $Z_0$.     

To prove that $\pi''$ is also smooth and equipped with local sections, we
consider the following extension $\pi''_0$ of $\pi''$ to the variety $X_0$ of
all pairs $(v,\psi) \in \fac \times K^{d \times r}$ satisfying $\psi \alpha(p)
= \alpha(v) \psi$ for $\alpha \in \Gamma^*$; in other words, the second
components of the pairs we single out run through the $\la$-homomorphisms $P
\rightarrow P/\rho(v)$.  Let $\pi'_0$ be the projection $X_0 \rightarrow \fac$
onto the first component.  Since X is contained in $X_0$ as an open subset and
$\pi'$ is the restriction of $\pi'_0$ to $X$, we need only establish our claim
for $\pi'_0$, and to do so, it suffices to check that the latter map is a
vector bundle projection.  But this follows again from \cite{\BongTrond, Lemma
3.4}, in view of the following two observations: (a) the set of
$\la$-homomorphisms $P
\rightarrow P /\rho(v)$ inside
$K^{d \times r}$ arises as the solution set of a homogeneous system of linear
equations, the coefficient matrix $A(v)$ of which depends morphically on
$v$  --  recall that $P$ is fixed; and (b), the  rank 
of $A(v)$ is locally constant on $\fac$.  Indeed the nullity of $A(v)$ equals
$\dim_K \Hom_\la (P, P/\rho(v))$.  Thus local constancy
can be gleaned from the fact that $\dim_K \Hom_\la(P,-)$ is
constant on the intersections of
$\fac$ with the open connected components of $\modla{d}$, because each such
component consists of modules having the same class in the Grothendieck group
${\roman K}_0(\lamod)$. $\qed$   
\enddemo

\proclaim{Proposition C}  There is an inclusion-preserving one-to-one
correspondence between the
$\aut_{\la} (P)$-stable subsets of $\sub$ and the $\GL_d$-stable subsets of
$\fac$ given by $M \mapsto \sigma\rho^{-1}(M)$.  Both ways, this correspondence
preserves openness, closures, connected and irreducible components, as well as
types of singularities.  Furthermore, an
$\aut_{\la} (P)$-stable subvariety $M$ of $\sub$ admits a geometric quotient
by $\aut_{\la} (P)$ if and only if the corresponding
$\GL_d$-stable subvariety of $\fac$ admits a geometric quotient by the
$\GL_d$-action; in the positive case, both of these quotients are
isomorphic to the geometric quotient of $\rho^{-1}(M)$ by $\aut_{\la} P
\times \H$. \endproclaim

\demo{Proof}  All of the statements except for the last are straighforward
consequences of Lemmas 1,2 and the general background given in
\cite{\BongTrond, Section 5}.  For the final assertion, we use the fact that
$\rho$ and
$\sigma$ are flat morphisms and apply
\cite{\BongTrond, Lemma 5.9}. $\qed$   \enddemo

We require one more ingredient for the proof of Theorem A.
Recall that, given any closed subgroup $H$ of an algebraic group $G$, the set $G/H$
of left cosets of $H$ in $G$ can be equipped with a structure of
quasiprojective variety  which makes the canonical surjection $G \rightarrow
G/H$ a geometric quotient of $G$ relative to the action of $H$ by right
translation (see, e.g., \cite{\Bor, Theorem 6.8}). Note that
$G/H$ carries a canonical transitive left $G$-action.   We will call a $G$-space $X$ a
{\it homogeneous\/}
$G$-space in case there exists a $G$-equivariant isomorphism $X \rightarrow
G/H$ for a suitable closed subgroup $H$ of
$G$.  As is well-known (see \cite {\Bor, Prop. 6.7}), provided that the
action of $G$ on $X$ is transitive, the orbit map
$G \rightarrow X$, $g \mapsto g x_0$ for $x_0 \in X$, induces an isomorphism
$G/\stab(x_0)
\rightarrow X$ if and only if this orbit map is separable.  In case
$K$ has characteristic zero, separability is automatic.    

The following structure theorem for homogeneous spaces of unipotent groups
is a consequence of a theorem of Rosenlicht \cite{\Ros, Theorem 1}.

\proclaim{Theorem (Rosenlicht)} Let $G$ be a unipotent
algebraic group. Then any homogeneous $G$-space is isomorphic to a full affine
space $\AA^r$. $\qed$
\endproclaim

\demo{Proof of Theorem A}  We postpone our argument for part (1) until
the end of the proof. 

Part (2). The second assertion under
(2) is well-known to follow from the first: Use
\cite{\Kra, II.2.6} to deduce the generic behavior of the fibre dimension.  To
obtain closedness of the $G$-orbits, we need only show that their intersections
with the patches $\guni(p)$ of our affine cover of $\gunis$ are all closed.  But
since each of the $\guni(p)$ is stabilized by $G$  --  see the discussion
preceding Theorem A  --  this is guaranteed by Kostant's result that the
orbits of a unipotent group acting on an {\it affine} variety are closed (c.f.
\cite{\Hum, p.115, Exercise 8}).     
   
Part (3).  Fix $C \in \gunis$.  In view of Rosenlicht's theorem, the claimed 
structure of the fibre $\F = G.C$ will follow if we can show $\F$ to be a
homogeneous
$G$-space.  It thus suffices to verify that the orbit map $G \rightarrow
\F$, $g\mapsto gC$, is separable.  To this end, it is enough to check
that the fibre $G_1$ of the orbit map is reduced (combine
\cite{\Bor, Prop. 6.7} with
\cite{\Kra, AI.5.5, Satz 2} to see this).  Now reducedness of
the stabilizer subgroup $G_1$ of $C$ is in turn automatic if $G_1$
arises as the solution set of a system of linear equations over $K$
(for background on reduced fibres  see \cite{\Kra, AI.2.5, 2.6}).  In our
situation, the latter condition can be verified as follows. 

Suppose that the uniserial module $\la e/ C$ has mast $p$, and let $p_0,
\dots, p_l$ be the right subpaths as before.  We know that there is a unique
family of scalars
$k_i(\alpha,u)$ such that $\alpha u - \sum_{i \in I(\alpha,u)}
k_i(\alpha,u) p_i \in C$.  Since, as a $\la$-module, $C$ is generated by these
differences and the non-routes on $p$, the condition $gC = C$ is clearly
equivalent to the requirement that
\smallskip 

\noindent $(\dagger)$ \qquad   all of the $g$-shifts $\ g\cdot (\alpha u -
\sum_{i
\in I(\alpha,u)} k_i(\alpha,u) p_i)\ $ for $(\alpha,u) \detour p$ belong to
$C$. 
\smallskip 

\noindent 
Let $B$ be a $K$-basis for $eJe$, i.e., 
$G^{\text{op}} = \{ e + \sum_{b \in B} l_b b \mid l_b \in K \}$. 
Using the fact that the images of the $p_i$ form a basis for 
$\la e / C$, we expand the elements $\alpha u b$ and $p_i b$ for $b \in B$ modulo
$C$.  A comparison of coefficients of the $p_i + C$ then shows condition
$(\dagger)$ to be tantamount to a linear system of equations for the
coefficients $l_b$ of $g$.   

To compute the dimension of the fibre $G.C$, observe that we have a group
epimorphism $\chi:G_1
\rightarrow \aut_{\la}^u(\la e/ C)$, where $\chi_g(\lambda + C) = g\lambda
 + C$.  Denoting by $G_0 \vartriangleleft G_1$ the kernel of $\chi$, we
thus obtain $\dim G_1 / G_0 = \dim \aut_{\la}(\la / C) - 1 = \dim_K \End_\la(U)
-1$.  It thus suffices to prove that $\dim G / G_0 = \mu(\SS) - 1$.  

Abbreviating
this last difference by $t$, we find precisely $t$ distinct oriented
cycles of positive length among the right subpaths $p_1, \dots, p_l$
of $p$; let these be $w_1, \dots, w_t$.  Moreover, let $x$ be the residue
class $e + C$, and consider the following set of top elements of $\la e /
C$, namely 
$$\T = \{(e + \sum_{i=1}^t c_i w_i) +C  \mid  (c_1, \dots,
c_t) \in \AA^t \}.$$
Clearly, $G$ acts transitively on $\T$, and
$G_0$ is precisely the stabilizer subgroup of $e+C$.  This gives us
$\dim G/G_0 = t$ as desired and completes the proof of the theorem.       

Part (1) We
will apply Proposition C to the situation $d=l+1$, the projective module $P = \la
e$, where
$e$ is the vertex going with the top simple of $\SS$, and the subvariety $M =
\gunis \subseteq \sub$.  Under the bijection
of Proposition C, the variety $M$ is paired with the $\GL_d$-stable
open subvariety $\moduni(\SS)$ of $\fac$, and consequently the first two claims under (1)
follow immediately from the proposition.  Concerning the final assertion, we
just need to know when $\aut_{\la}(\la e)$ acts trivially on
$\gunis$. But, keeping in mind that the orbits of the 
$\aut_{\la}(\la e)$-action coincide with those of the $G$-action, we obtain
the answer from part (3) of Theorem A:  This happens if and only if
$\mu(\SS) - \dim_K \End_{\la} U = 0$ for all uniserial modules $U$ with
composition series $\SS$. $\qed$ \enddemo 

\noindent {\bf Remarks:}
 
{\bf 1.} In view of Theorem A, Proposition C tells us
that the $\GL_d$-orbits of $\moduni(\SS)$ are closed in $\moduni(\SS)$, in
other words that uniserial modules have no proper uniserial
degenerations, and we thus re-encounter Proposition E of
\cite{\GeomII}.

{\bf 2.} As we pointed out earlier, our proof for the first part of Theorem A is
derived from a geometric re-interpretation of an obvious module-theoretic
correspondence: namely, that between the $\aut_{\la}(P)$-orbits of submodules of
a fixed dimension $m$ inside a projective module $P$ and the isomorphism types
of factor modules of dimension
$\ \dim P - m\ $ of
$P$.  There are other time-honored ways of studying modules which
translate into similar geometric pictures.  For instance, if one
focuses on the maps $P_1 \rightarrow P_2$ of fixed rank between
two projective modules $P_1$ and $P_2$, then the cokernels of two maps
$f$ and $g$ are isomorphic if and only if $f$ and $g$ are conjugate
under the obvious $\aut(P_1) \times \aut(P_2)$-action on
$\Hom(P_1,P_2)$.  To obtain the corresponding result at the geometric
level, one can again use some canonical bundle constructions modeling
our approach to Proposition C.  Analogous correspondences arise if one
interprets a module as the image of a map from a projective to an
injective object or as the kernel of a map between two injectives.  We
leave the details of the ensuing geometric setups to
the reader. 

{\bf 3.}  Implicitly, part of the relationship between Grassmannians and module
varieties developed above already played a role in
Gabriel's article `{\it Finite representation type is open}'; see
\cite{\Gab}.  To elaborate the connection a bit, we denote by
$\algs$ the affine variety of $s$-dimensional associative $K$-algebras
with identity, by $\algsmodr$ the variety of all pairs $(a,v)$ composed
of a point $a$ of $\algs$ and an $r$-dimensional
$a$-module $v$, and by $\algssubm$ the variety of pairs $(a,u)$, where $m =
(s-1)r$ and the second component is an $m$-dimensional submodule
$u$ of the $r$-th power of the regular representation of $a$; somewhat
sloppily we denote the latter by $a^r$.  Finally, $\bmr$ will stand for the
variety of pairs $(a,b)$, consisting of $a \in \algs$ and matrices $b$ running
through the elements of $\GL_{rs}$ with the property that
the first $m$ columns span a submodule of $a^r$.  This
setup comes with morphisms $\rho : \bmr \rightarrow \algssubm$ and
$\sigma : \bmr \rightarrow \algsmodr$ defined as before. 
They can be supplemented to a commutative diagram 

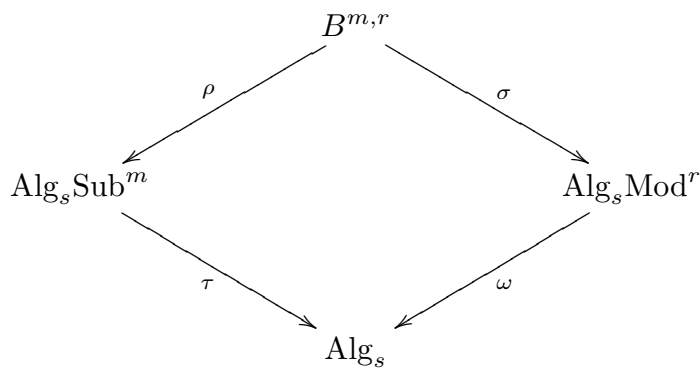
\ignore{
$$\xymatrixcolsep{4pc}\xymatrixrowsep{3pc}
\xymatrix{
 &\bmr  \ar[dl]_\rho \ar[dr]^\sigma \\
\algssubm  \ar[dr]_\tau &&\algsmodr \ar[dl]^\omega \\
 &\algs
}$$
}

\noindent where $\tau$ and $\omega$ are the obvious projections.  Again
$\rho$ and $\sigma$ are smooth morphisms, both
equivariant with respect to the appropriate group actions; 
moreover, $\tau$ is proper.  A crucial difference, compared
with the setup of Proposition C, however, lies in the fact
that the group $\aut(P)$ needs to be replaced by the obvious
smooth group scheme over $\algs$ consisting of automorphisms
of $a^r$.

One of the key points proved in Gabriel's paper is the fact that
$\omega$ maps any $\GL_{r}$-stable closed subset $X$ of
$\algsmodr$ to a closed subset $Y$ of $\algs$.  This can be
gleaned from our picture, because $Y = \tau \rho \sigma^{-1}
(X)$:  Indeed, $\rho$ is a principal $\H$-bundle, where $\H$ denotes the
subgroup of $\GL_{rs}$ relevant to the present situation, and consequently
$\rho$ maps the closed $\H$-stable subset
$\sigma^{-1}(X)$ of $\algsmodr$ to a closed subset of
$\algssubm$; the proper morphism $\tau$ then takes $\rho
\sigma^{-1}(X)$ in turn to a closed set.  Gabriel's original argument is
essentially this, but does not specify the `intermediate' variety
$\bmr$.  However, Kraft introduced an auxiliary
variety similar to ours in his treatment of Gabriel's
theorem in \cite{\Kraft}.

\head 3. Applications of Theorem A\endhead

\subhead{(I) Computational Remarks and an example}\endsubhead
We begin with remarks indicating how the quantities appearing in Theorem
A, such as the dimension of the fibre above a class $\Phi_p(k)$ of uniserial modules, 
can be determined algorithmically from quiver and relations of
$\la$.  Subsequently, we will give an example that serves a twofold purpose: 
Namely, (a) it shows that the fibre dimension can assume multiple values on a given
irreducible component of
$\gunis$, and (b) it displays constraints which, in general, render it
impossible to translate the described $G$-action on the
variety $\guni(\SS)$ into a more `traditional' conjugation action of a
unipotent matrix group on the affine variety $\uni(\SS) = \bigcup\uni(p)
\subseteq \modla{l+1}$, where $p$ traces the paths through $\SS$.  Such a
translation is not even possible for the individual affine patches $\guni(p)$ and
$\uni(p)$, respectively (keep in mind that they are unions of orbits, so that
our $G$-action on $\gunis$ restricts to actions on the $\guni(p)$).  More
precisely, while according to  Section 1,
$\guni(p)$ and the closed subvariety $\uni(p)$ of $\modla{l+1}$ are
always isomorphic, in general there is no isomorphism
$\guni(p)
\rightarrow \uni(p)$ which reduces the $G$-action on $\guni(p)$, via
some group homomorphism $G \rightarrow \GL_{l+1}$, to the conjugation
action of a unipotent subgroup $H \le \GL_{l+1}$ on $\uni(p)$. 
		
{\bf 1.} First we give a compact review of a method to determine
defining polynomials for the affine varieties
$V_p \cong \guni(p)$ from quiver and relations of $\la$; for more detail we
refer to \cite{\GeomI}. 

Let $\la = K\Gamma/I$ and $p$ a path in $\Gamma$.  Consider
the noncommutative polynomial ring $K\Gamma [X_i(\alpha,u) \mid (\alpha,u)
\detour p,\ i \in I(\alpha,u)]$; here, the variables
$X_i(\alpha,u)$ commute with the coefficients in $K\Gamma$ and with
each other.  We introduce a congruence relation 
relative to addition and left multiplication on this ring, denoted by
`$\seq$' and generated by the following {\it substitution equations}.  Namely,
the congruences 
$$\alpha u \seq  \sum_{i \in I(\alpha,u)} X_i(\alpha,u) v_i(\alpha,u) \qquad
\text{and} \qquad q
\seq 0$$   
for $(\alpha,u) \detour p$ and for those paths $q$ which fail to be
routes on $p$.  Given any element $z \in K\Gamma$, there are unique elements
$\tau_i(z)$ in the commutative polynomial ring $K[X_i(\alpha,u) \mid (\alpha,u)
\detour p,\ i \in I(\alpha,u)]$ with the property that 
$$z \seq \sum_{i=0}^l \tau_i(z) p_i.$$ 
Provided one starts with a set of
generators for the left ideal $Ie(0)$ of
$K\Gamma$, say  $z_1, \dots z_r$, the corresponding polynomials
$\tau_i(z_j)$, $0 \le i \le l$, $1 \le j \le r$, determine $V_p$;  namely, $V_p =
V(\tau_i(z_j) \mid i,j)$.  In practice, these polynomials are obtained by
starting with a choice of $z_j$ and successively inserting the substitution
equations into monomials of the form $f(X) v$, where $f(X) \in K[X]$ and $v \in
K\Gamma$ is a path (for more detail, see \cite{\GeomI}).
 
{\bf 2.} Next we show how to obtain the fibre dimensions of the maps
$\Phi_p$  --  or, equivalently, of the maps $\phi_{\SS}$. In view of Theorem
A, this boils down to determining the rank of a certain $t\times
t$-matrix $A(k)$, where $t = \mu(\SS) - 1$. 

Given any point
$k\in V_p$, the dimension
$\delta(k) = \dim_K\End_{\la}(U)$, where $U$ is any representative of the
isomorphism class $\Phi_p(k)$, can be computed along the following line.
(Actually, our instructions can be used to explicitly calculate
$\aut^u_\la(U)$.)  
Suppose $p$ passes through the sequence $\SS$ with top simple $S(0) = \la e/
Je$, and let $w_0 = e, w_1 \dots, w_t$ be the different right
subpaths of $p$ ending in the starting vertex $e$; clearly, $\mu(\SS) =
t+1$ by the definition of $\mu$.  For a  suitably chosen top element
$x = ex$ of $U$, we then have
$$\alpha ux= \sum_{i\in I(\alpha,u)} k_i(\alpha,u) v_i(\alpha,u) x
\tag \dag$$  for all $(\alpha,u)\detour p$. We indicate how
to set up a system of linear equations for those $t$-tuples
$(c_1,\dots,c_t)$ of scalars for which the assignment
$$x\mapsto x + \sum_{i=1}^t c_iw_ix$$ induces a unipotent
$\la$-automorphism of $U$: We repeatedly use the equations
$(\dagger)$ to expand the elements $\alpha uw_jx$ and
$v_i(\alpha,u)w_jx$, for $(\alpha,u)\detour p$, $i\in I(\alpha,u)$ and
$j=1,\dots,t$, in terms of the $K$-basis $x=p_0x, p_1x,
\dots, p_lx$ of
$U$. Next we insert these expansions into the equations
$$\alpha u \biggl( x+ \sum_{j=1}^t c_jw_jx \biggr) =\sum_{i\in
I(\alpha,u)} k_i(\alpha,u) v_i(\alpha,u) \biggl( x+ \sum_{j=1}^t c_jw_jx
\biggr) \tag \ddag$$ for $(\alpha,u)\detour p$, and compare coefficients of
the $p_jx$'s so as to obtain a system of linear equations for
$c_1,\dots,c_t$. This system is homogeneous, say of the form
$A(k)(c_1,\dots,c_t)^T\allowmathbreak =0$ for a suitable
$K$-matrix $A(k)$, as the `constant' terms on the two sides of
$(\ddagger)$ obviously cancel out. Now $\delta(k)= t-\rank A(k)$,
and consequently the dimension of the fibre $\Phi_p^{-1}\Phi_p(k)$ equals the
rank of the matrix $A(k)$.

\subhead Example 1 \endsubhead Fixing $l \in \NN$, let $\la =
K\Gamma/I$ be the algebra based on the quiver $\Gamma$

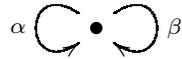
\ignore{
$$\xymatrixcolsep{4pc}
\xymatrix{
\bullet \ar@(ul,dl)_\alpha \ar@(ur,dr)^\beta
 }$$ 
}

\noindent and the ideal $I \subseteq K\Gamma$ generated by all paths of length
$l+1$.  Moreover, consider the $(l+1)$-term sequence $\SS = (S,S, \dots,S)$,
where
$S$ is the unique simple left $\la$-module.  Given any path $p$ of length $l$,
the first of the above remarks
permits us to compute $\guni(p) \cong V_p \cong \AA^{l(l+1)/2}$.  According to
\cite{\GeomII, part (3) of Example}, we moreover have $\gunis \cong
\PP^1 \times \cdots \times \PP^1 \times \AA^{l(l+1)/2}$ with $l$ factors $\PP^1$
occurring in the product; in particular,
$\gunis$ is irreducible.

First we will see that all of the numbers between $0$ and $l-1$ show up as fiber
dimensions.  For that purpose consider, for instance, the path $p =
\alpha^{l-1}\beta$ through $\SS$, and for $1 \le m \le l$, let $k(m) \in V_p$ be
a point corresponding to the uniserial module
$$U(m) = \la /\bigl( \la(\alpha - \beta^{m-1}\alpha)\ + \ \sum_{0 \le j \le l-1}
\la\beta \alpha^{j} \beta\bigr).$$ 
In other words, the detours on $p$ being $(\alpha, 1)$ and $(\beta,
\alpha^j\beta)$ for $j \ge 1$, the point 
$k(m)$ may be chosen to have coordinates
$k_m(\alpha,1) = 1$, $k_i(\alpha,1) = 0$ for $i \ne m$, and
$k_i(\beta,\alpha^j \beta) = 0$ whenever $j \ge 0$.  Using Theorem A
and Remark 2 above, it is now easy to
compute the dimension of the fibre of $\Phi_p$ (or, equivalently, that of the
fibre of $\phi_{\SS}$) over $U(m) = \Phi_p(k(m))$ to be $l-m$. 

Next we show that there is no unipotent subgroup $H$ of $\GL_{l+1}$ which
acts nontrivially on $\uni(p)$ by conjugation.  This of
course entails that there is no conjugation action of a unipotent matrix group 
$H$ on $\uni(p)$ such that the orbits of the action coincide with the fibres of
the canonical map from
$\uni(p)$ to the isomorphism classes of uniserials with mast $p$.  To verify our claim,
observe that, in our present setting, $\uni(p)$ consists of all pairs
$(A,B)$ of strictly lower triangular $(l+1)\times(l+1)$-matrices representing
multiplication by $\alpha$ and $\beta$, with the property that the first column
of
$B$ consists of the canonical basis vector $e_2$, and the second through
$l$-th columns of
$A$ coincide with $e_3$ through $e_{l+1}$, in that order.  Now it is
readily checked that the only nontrivial unipotent subgroup of
$\GL_{l+1}$ which leaves $\uni(p)$ invariant under conjugation is $\{E +
cE_{l+1,1} \mid c \in K\}$;  here $E$ is the identity matrix and the $E_{ij}$
denote the matrix units in $\text {M}_{l+1}(K)$.  Since
$\uni(p)$ is made up of strictly lower triangular matrices, this group clearly
acts trivially on $\uni(p)$. $\qed$

\subhead{(II)  Algebras of finite uniserial type revisited}\endsubhead
In \cite{\Bong}, the first author settled the following conjecture from
\cite{\GeomIII} affirmatively:  If
$\la$ has finite uniserial type, then, given any sequence $\SS$ of simple
$\la$-modules, there is at most one uniserial module (up to isomorphism)
having composition series $\SS$. We start by
indicating a short geometric proof of this fact and then combine it with
Theorem A to confirm the remaining conjecture of
\cite{\GeomI} as well:  Namely, finite uniserial type of $\la$ forces each
nonempty variety
$V_{\SS}$ to be a full affine space.  Our main application,
Corollary D, will provide new equivalent characterizations of such
algebras which round off the pictures presented in \cite{\GeomIII} and \cite{\Bong}. 
All the
quantities involved in these descriptions can be computed from quiver and
relations of $\la$. 
 
As before, let $\SS = (S(0), \dots,
S(l))$ be a sequence of simple left
$\Lambda$-modules and set $d = l+1$.  Recall that the set $\moduni(\SS)
\subseteq \modla{d}$, consisting of the points of $\modla{d}$ that correspond to
uniserial modules $U$ with composition series $\SS$, is an open subvariety of
$\modla{d}$ which is stable under the
$\GL_d$-action.   

\proclaim{Bongartz's Theorem} {\rm \cite{\Bong, Theorem 1}} If $\la$ has finite
uniserial type, then, given any finite sequence $\SS$ of simple left
$\la$-modules, there exists at most one uniserial left $\la$-module with
composition series
$\SS$, up to isomorphism. \endproclaim

\demo{Proof} Assume finite uniserial type.  Keeping the above notation, we
perform an induction on $l$. Clearly, the case
$l=0$ is harmless, so suppose that $l\ge1$ and that $\gunis$ is
nonempty; the latter is tantamount to the requirement that the
quasi-affine variety
$\moduni(\SS)$ be nonempty.  By induction hypothesis, all of the uniserial
modules
$U$ with composition series
$\SS$ have the same radical, up to isomorphism  --  call it $R$. Consequently,
all uniserial modules with composition series $\SS$ are extensions of $R$ by
$S(0)$. Now it is known that the points in $\modla{d}$ corresponding to
arbitrary extensions of $R$ by $S(0)$ form an irreducible subset
$E(R,S(0))$ of $\modla{d}$ (see, e.g., \cite{\BongAdv, 6.3}), and since
$\moduni(\SS)$ is open, the same is true for the intersection $E(R, S(0)) \cap
\moduni(\SS)$. Investing the fact that
$\la$ has finite uniserial type, we infer that this intersection is a
{\it finite} union of $\GL_d$-orbits. All of these orbits being closed in
$\moduni(\SS)$ by Remark 1 at the end of Section 2, irreducibility implies that there is
only one such orbit, which is what we wanted to show.
\qed\enddemo

Again suppose that the primitive idempotent corresponding to the top simple
$S(0)$ of $\SS$ is $e$ and that $p$ is a path through $\SS$. For any point
$C \in\gunis$ and any
$k \in V_p$, we denote by $\delta(C)$ and $\delta(k)$ the $K$-dimensions of the
$\la$-endomorphism rings of the uniserial modules $\la e / C$ and
$\Phi_p(k)$, respectively. Recall that $\mu(\SS)$ is the
multiplicity of $S(0)$ in $\SS$; accordingly, we write $\mu(p)$ for $\mu(\SS)$ whenever
$p$ is a path through $\SS$. We then obtain the following characterizations of finite
uniserial type.  On the side, we note that, since polynomials for the varieties $V_p$ can
be directly obtained from quiver and relations of a split basic algebra, condition (3) is
verifiable by means of Gr\"obner basis methods.  The formally strongest condition, (4), 
demonstrates how restrictive the condition of finite uniserial type is.       

\proclaim{Corollary D to Theorem A} The following statements are equivalent for
$\la= K\Gamma/I$:

{\rm (1)} $\la$ has finite uniserial type.

{\rm (2)} For each finite sequence
$\SS$ of simple left
$\la$-modules and \underbar{every} point $C\in
\gunis$, the variety $\gunis$ is either empty or has dimension 
$\mu(\SS) - \delta(C)$.

{\rm (3)} $\Gamma$ has no double arrows, and for each path $p \in K\Gamma$,
the affine variety
$V_p$ is either empty or irreducible of dimension $\mu(p) -\delta(k)$ for
\underbar{some} point
$k\in V_p$.

{\rm (4)} If $\SS$ is the sequence of consecutive composition factors of a
uniserial left $\la$-module, then there exists a unique path $p$ through $\SS$
and $\gunis = \guni(p) \cong V_p$ is isomorphic to
$\AA^{\mu(\SS) -\delta(k)}$ for all points $k \in V_p$.  

\endproclaim

Of course, all statements of the corollary are equivalent to their
right-hand analogues.

\demo{Proof} (1)$\implies$(2) is a consequence of Bongartz's theorem, reproved above,
and Theorem A. For (2)$\implies$(1), suppose that $p$ is a path through
$\SS$ with $\guni(p) \ne \varnothing$ and $\C$ an irreducible component of
$\guni(p)$, containing a point $C$ say.  Then $\dim \C \le \mu(\SS) - \delta(C)$
by (2), and by Theorem A, the fibre of $\Phi_\SS$ over the isomorphism class of
$\la e/C$ is a closed subvariety of $\C$.  Moreover, Theorem A shows this fibre
to have dimension 
$\mu(\SS) - \delta(C)$, making it coincide with $\C$.  This establishes
finiteness of the number of isomorphism classes of simples with composition
series $\SS$ as required.   

(1)$\implies$(4).  Clearly, finite uniserial type excludes the existence
of double arrows in $\Gamma$, so that, whenever $\gunis \ne \varnothing$,
there exists a unique path $p$ through $\SS$, and $\gunis = \guni(p)
\cong V_p$.  The rest follows again from Bongartz's theorem and Theorem A.

The implication (4)$\implies$(3) is obvious, 
and the reason for (3)$\implies$(1) is analogous to that given for
(2)$\implies$(1); note that the irreducibility condition in (3)  permits us to
cut the dimension hypothesis down to a single point $k \in V_p$.   
\qed\enddemo
\medskip

\noindent {\it Applying Corollary D to Examples}

In case one wishes to decide whether a given algebra $\la= K\Gamma/I$ has
finite uniserial type, it is most efficient to start by checking the
following strong necessary condition (N) given in \cite{\GeomIII, Section 3}: 

(N) Whenever
$\alpha : e\rightarrow e'$ is an arrow in $\Gamma$ and $p : e\rightarrow e'$
a mast of positive length, $p$ is of the form

\ignore{
$$\xymatrix{ e \ar[rr]^\alpha &&e' \udotloopr{c'} &&\text{or}&& e
\ar[rr]^\alpha
\udotloopr{c} &&e' }$$ }

\noindent where $c',c$ are oriented cycles which may be trivial.  
Given that condition (N) is satisfied, the results of
\cite{\GeomIII, Sections 4,5} permit us in many cases to decide `at a glance'
whether, for a given path $p$, there are infinitely many uniserials with
mast $p$, up to isomorphism. For the remaining paths, the decision process is
most effortlessly carried out with the help of conditions (2), (3) of Corollary D, combined
with the concluding remarks of Section 2.  We illustrate the procedure with two algebras. 
For these, answering the question `finite or infinite uniserial type?', based on previous
methods, is cumbersome.

\subhead Examples 2\endsubhead 
Consider the following quiver $\Gamma$: 
\ignore{
$$\xymatrixcolsep{3pc}
\xymatrix{ 1 \ar[r]<0.75ex>^{\alpha_1} &2 \ar[l]<0.75ex>^{\alpha_2}
\ar[r]<0.75ex>^{\alpha_3} \ar[d]^{\alpha_5} &3
\ar[l]<0.75ex>^{\alpha_4}\\
 &4
 }$$ }
\noindent Clearly, each finite dimensional factor algebra $\la$ of $K\Gamma$
satisfies condition (N). In particular, absence of double arrows permits us
tho identify the variety $\gunis$ with $V_p$, whenever $p$ is a path through
a sequence $\SS$ of simple $\la$-modules.

(a) Let $\la= K\Gamma/I$, where $I$ is generated by
$\alpha_2\alpha_1\alpha_2$,
$\alpha_3\alpha_4$, $\alpha_2\alpha_4$, and $\alpha_5\alpha_1\alpha_2
-\alpha_5\alpha_4\alpha_3\alpha_1\alpha_2$. Then it is easy to see that the
uniserial dimension of $\la$ at $\SS'$ is $\le 0$ for all $\SS'\ne \SS=
(S_1,S_2,S_1,S_2,S_3,S_2,S_4)$. We will use the computational
remarks at the end of Section 2 to check that $\uniserdim=1$, which, in
particular, means that
$\la$ fails to have finite uniserial type.

The unique path $p$ of length 6 through $(e_1,e_2,e_1,e_2,e_3,e_2,e_4)$ is
$p=
\alpha_5\alpha_4\alpha_3\alpha_1\alpha_2\alpha_1$. To compute $V_p$, we note
that the detours on $p$ are
$(\alpha_5,\alpha_1)$,
$(\alpha_3,\alpha_1)$, and
$(\alpha_5,\alpha_1\alpha_2\alpha_1)$; they lead to the substitution
equations
$$\alpha_3\alpha_1 \seq X_1\alpha_3\alpha_1\alpha_2\alpha_1, \qquad
\alpha_5\alpha_1 \seq X_2p, \qquad \alpha_5\alpha_1\alpha_2\alpha_1 \seq
X_3p.$$

The relation $\alpha_5\alpha_1\alpha_2
-\alpha_5\alpha_4\alpha_3\alpha_1\alpha_2$ can be supplemented by non-routes
on $p$ to a generating set for the left ideal $Ie_1$, whence insertion of
the substitution equations yields $X_3=1$ and no conditions impinging on
$X_1,X_2$. Thus, $V_p= \AA^2$. The matrix $A(k)$ of Remark 2 of Section 2 has
rank 1 for each
$k= (k_1,k_2)\in V_p$, whence the fibre dimension of
$\Phi_p^{-1}\Phi_p(k)$ is constant on $V_p$ and equal to 1. Thus
the uniserial dimension of $\la$ at $\SS$ is indeed $1$. 
\medskip

(b) Now let $\la'= K\Gamma/I'$, where $I'$ is generated by
$\alpha_2\alpha_1\alpha_2$,
$\alpha_3\alpha_4$, $\alpha_2\alpha_4$, and $\alpha_5\alpha_1
-\alpha_5\alpha_4\alpha_3\alpha_1$, and let $\SS$ and $p$ be as in part (a).
This time, the relations $\alpha_5\alpha_1\alpha_2\alpha_1 -p$ and
$\alpha_5\alpha_1 -\alpha_5\alpha_4\alpha_3\alpha_1$ can be supplemented by
non-routes on $p$ to a generating set for the left ideal $I'e_1$, and our
substitution equations yield $X_3=1$ and $X_1=X_2$. Thus, $V_p\cong \AA^1$.
Again we compute $\rank A(k)=1$ for arbitrary $k\in V_p$, which implies that
the uniserial dimension of $\la'$ at $\SS$ equals 0 in this case.  The
other $V_q$'s can be dealt with in a computation-free manner (see
\cite{\GeomIII}), the conclusion being that $\la'$ does have
finite uniserial type.
\qed

\head 4. Endomorphism rings of uniserial modules and an alternate road to the fibre
structure
\endhead

In the following, we will present an alternate approach to the fact that the
fibres of the maps
$\Phi_p$ are full affine spaces, thereby encountering the somewhat surprising
fact that all endomorphism rings of uniserial
$\la$-modules are commutative (keep in mind that our blanket hypothesis on
$\la$ calls for an algebraically closed base field).

Again $p = \alpha_l
\cdots \alpha_1$ will be a fixed path of length $l$ passing through
the sequence $(e(0), \dots, e(l))$ of vertices that accompanies our  sequence
$\SS = (S(0), \dots, S(l))$ of simples.  Moreover, we will continue to
identify the preferred primitive idempotents $e_1, \dots, e_n \in \la$ with
the vertices of $\Gamma$.  This time we will work with the isomorphic
copy $\uni(p)$ of
$\guni(p)$, alias
$V_p$, inside the classical variety $\moduni(\SS)$, as introduced at
the end of Section 1.  Namely,
$$\uni(p) = \{x = (\alpha(x))_{\alpha\in \Gamma^*} \in \moduni(\SS) \mid
\alpha_i(x) \bold e_{i-1} =
\bold e_i \ \text{for}\ 1 \le i \le l \},$$ 
where $(\bold e_0, \dots \bold e_l)$ is the canonical basis for $K^{l+1}$, and
$\Gamma^*$ is the union of
$\{e_1,\dots,e_n\}$ with the set of arrows of $\Gamma$.  Here the
components
$\alpha(x)$ of a point $x$ are viewed as $(l+1) \times
(l+1)$-matrices representing multiplication by the pertinent $\alpha$ relative
to the canonical basis;  the rows and columns are indexed by
$\{0,
\dots, l\}$.  In light of the commutative diagram given at the end of Section 1,
which links up all of our canonical maps from varieties to module categories, it is
harmless to identify $V_p$ with the subvariety
$\uni(p)$ of $\moduni(\SS)$, and $\Phi_p$ with the restriction $R|_{\uni(p)}$   -- 
we will do this in the sequel.  Note that, for any
$x \in V_p$, there is a uniserial module representing $\Phi_p(x)$ having
underlying vector space
$K^{l+1}$ and composition series $\la \bold e_0 \supset \la \bold e_1
\supset \dots \supset \la \bold e_l \supset 0$; this is a trivial
consequence of the definitions.  Hence all of the matrices
$\alpha(x)$,  where $\alpha$ runs through the arrows of $\Gamma$, are
strictly lower triangular.

We will denote by $\boldd$ the sequence $(d_1,\dots,d_n)$, where $d_i$ is the
multiplicity of the simple module $\la e_i/Je_i$ in $\SS$, and by $X_i$ the
$d_i$-dimensional subspace of $X= K^{l+1}$ spanned by those $\bold e_j$ for
which $S(j)\cong \la e_i/Je_i$. Then, clearly, $X= \bigoplus_{1\le i\le n}
X_i$, and $\moduni(\SS)$ is contained in the closed subset $\modla{\boldd}$ of
$\modla{l+1}$ which consist of the points $(\alpha(x))_{\alpha\in \Gamma^*}$
with $\alpha(x)\in \operatorname{Hom}_K (X_{s(\alpha)}, X_{t(\alpha)})$,
where $s(\alpha)$ and $t(\alpha)$ denote the starting and terminal vertices
of $\alpha$, respectively; by this we mean in particular that the restriction of
$\alpha(x)$ to $\sum_{i \ne s(\alpha)} X_i$ is zero.  Clearly,
$\modla{\boldd}$ is stable under conjugation by $\GL(\boldd)= \prod_{1\le i\le
n} \GL_{d_i}$ if the elements in $\prod_{1\le i\le n} \GL_{d_i}$ represent those
automorphisms of
$X$ which leave the
$X_i$ invariant, with respect to the canonical basis $\bold e_0, \dots, \bold e_l$ of $X$.
Note that the image of $\modla{\boldd}$ under $R$ is just the
set of isomorphism classes of $\la$-modules with dimension vector $\boldd$.

From now on, we will assume that the starting point $e(0)$ of $p$ equals
$e_1$ in the above listing of the vertices. This entails, in particular, that the first
canonical basis vector
$\bold e_0$ of $X$ coincides with the first canonical basis vector of $X_1$.

 We will, moreover, consider two subgroups of the linear group
$\GL(\boldd)$, the first being $U(\boldd) = \prod_{1 \le i \le
n} U(d_i)$, where
$U(d_i)$ consists of the unipotent lower triangular matrices in $\GL(d_i)$. 
While the group
$U(\boldd)$ acts on $\modla{\boldd}$ by conjugation, it also acts 
on the vector space $X$ by multiplication from the left.  The
identification of elements of $U(\boldd)$ with $(l+1) \times (l+1)$-matrices
is to follow the same pattern as above. Secondly, we consider the stabilizer
subgroup
$H$ of $\bold e_0$ in $U(\boldd)$ under the latter action; i.e., $H$ is the
subgroup of those matrices in $U(\boldd)$ which carry $\bold e_0$ in the $0$-th
column.

Observe that the variety $V_p \subseteq \modla{\boldd}$ is not stable under
conjugation  by elements in $U(\boldd)$.  But the following subvariety $V_p'$ of
$\modla{\boldd}$ containing
$V_p$ is: Namely, let $V_p'$ consist of those points $x$, represented by
tuples of lower triangular matrices $\alpha(x) \in M_{l+1}(K)$, for which
$\alpha_i(x)_{i,i-1} = 1$.  In particular, $V_p'$ is stable relative to
conjugation by elements in $H$.  

\proclaim{Lemma 3}  Restriction of the $H$-action to $V_p$ induces an
isomorphism
$$H \times V_p \rightarrow V_p'.$$ \endproclaim

\demo{Proof} The inverse of the morphism  $(h,x) \mapsto h.x$ is given by  
 $$y \mapsto (g(y),g(y)^{-1}yg(y)),$$ where the columns of $g(y)$ are
$\bold e_{0},\alpha_{1}(y) \bold e_{0},
\ldots,\alpha_{l}(y)\alpha_{l-1}(y) \ldots \alpha_{1}(y) \bold e_0$.
\qed\enddemo  

We now fix an element $x \in V_p$.  As before, we let $w_1, \dots, w_t$ be the
distinct right subpaths of positive length ending in the starting vertex
$e$ of $p$, insisting that they be listed in order of
increasing length.  We write them as successive extensions of one another:  Set
$u_1 = w_1$, and let $u_i$ be the path with
$w_i$ = $u_i w_{i-1}$.  Note that we
have $\dim X_1 = d_1 = t+1$.  Clearly the maps $u_j(x)$ induce nilpotent
endomorphisms of $X_1$, represented by strictly lower triangular matrices,
while sending $\sum_{2 \le i \le n} X_i$ to zero. The definition of $V_p$
moreover guarantees that
$w_t(x) = u_t(x) u_{t-1}(x) \dots u_1(x) \ne 0$ and  --  our base field $K$
being infinite  --  we can therefore find a $K$-linear combination $f$ of the
maps $u_1(x), \dots, u_t(x)$ such that $f^t$ does not vanish.  Let $A
\subseteq U(\boldd)$ be the subgroup consisting of those elements of $U(d_1)
\times 1 \times \dots
\times 1$ which commute with $f$.  Elementary linear algebra then shows that 

\noindent $\bullet$ \quad $A$ consists precisely of the linear combinations
of positive powers of
$f$ added to the identity matrix in $U(\boldd)$  -- so, in particular, $A$ is
commutative  --  and

\noindent $\bullet$ \quad $A \bold e_0 = X_1$ inside $K^{l+1} = X_1 \oplus
\dots \oplus X_n$.

\proclaim{Lemma 4} The varieties $H \times A$ and $U(\boldd)$  are isomorphic
via the multiplication map $(h,a) \mapsto ha$. \endproclaim

\demo{Proof}  We observe that $H \cap A = 1$, since all nonzero linear
combinations of $f, f^2, \dots, f^t$ are strictly lower diagonal with
nonzero first columns.  Combined with the second of the preceding
observations, this shows that, for each $g \in U(\boldd)$, there exists a unique
element $a(g) \in A$ with the property that $g \cdot a(g)$ belongs to
$H$, i. e., with $a(g) \bold e_0 = g^{-1} \bold e_0$.  Next we note that the
assignment $U(\boldd)
\rightarrow A$, $g \mapsto a(g)$ is a morphism of varieties;  indeed, $a(g)
= id +
\sum_{i=1}^t c_i f^i$, where the $c_i$ are polynomials in the entries of the
first column of
$g^{-1}$, with coefficients depending only on $f$.  Consequently, the
inverse to the multiplication map of our claim, namely the map $U(\boldd)
\rightarrow H \times A$, $g \mapsto (g \cdot a(g),a(g)^{-1})$, is a morphism
as well. \qed \enddemo

We are now in a position to prove the following modified version of Theorem
A.  

\proclaim{Theorem A$'$}  Let $U$ be a uniserial module representing the
isomorphism class
$\Phi_p(x)$. Then the following are true:
\medskip

\noindent {\rm(1)}  The group $\aut_{\la}^u(U)$ of unipotent
$\la$-automorphisms of $U$ is isomorphic to a subgroup $B$ of $A$. 
Consequently, all endomorphism algebras of uniserial
$\la$-modules are commutative. 
\medskip

\noindent {\rm(2)}  The fibre $\Phi_p^{-1}\Phi_p(x)$ is isomorphic to the
commutative unipotent group $A/B$ and, in particular,  this fibre is
isomorphic to $\AA^m$ with $$m = t - \dim
\aut_{\la}^u(U) = t+1 - \dim \End_{\la}(U).$$ 
\endproclaim

\demo{Proof}  We may assume that the module $U$ has underlying space $X =
K^{l+1}$ and composition series $\la \bold e_0 \supset \la \bold e_1
\supset \dots \supset \la \bold e_l \supset 0$.  Then we obtain
$e_iU = X_i$ for $1 \le i \le n$, and may view $\aut_{\la}^u(U)$ as a
subgroup of $U(\boldd)$.  

(1)  By definition, every element $g \in U(\boldd)$ can be uniquely written in
the form $(g_1,g_2)$, where $g_1 \in U(d_1)$ is a unipotent automorphism of
$X_1$ and $g_2 \in \prod_{i \ge 2} U(d_i)$ a unipotent automorphism of
$\sum_{i \ge 2} X_i$.  If $g \in \aut_{\la}^u(U)$, then $g$ commutes with
the $u_j(x)$ and hence with $f$.  This means that $(g_1,id)$ commutes with
$f$, i. e., that $(g_1,id) \in A$.  Since $\bold e_0$ is a top element of
$U$, each $\la$-automorphism
$g \in \aut_{\la}^u(U)$ is determined by $g(\bold e_0) = g_1(\bold e_0)$,
and we deduce that the map
$\aut_{\la}^u(U) \rightarrow A$, given by $g \mapsto (g_1,id)$ is a
monomorphism of groups.  For the final assertion of of part (1), we note that
each element 
$f \in \End_{\la}(U)$ is of the form $f = a\cdot id_U + g$, for some scalar $a$
and $g \in \aut_{\la}^u(U)$, and thus infer commutativity of $\End_{\la}(U)$
from that of $\aut_{\la}^u(U)$.

(2)  Let $B \subseteq A$ be the image of the group homomorphism under (1).  
We begin by observing that $\Phi_p^{-1}\Phi_p (x) = O(x) \cap V_p$, where
$O(x)$ is the
$U(\boldd)$-orbit of $x$ in $V_p'$.  We will show in two steps that this
intersection is a geometric quotient of $U(\boldd)$ and will subsequently
ascertain that it is indeed isomorphic to $A / B$.  Since $A /
B$ is a connected unipotent group, Theorem 1 will then guarantee that $O(x)
\cap V_p \cong A / B \cong \AA^m$, and our claim will follow from the fact
that $\dim A = t$.     

On one hand, $O(x)$ is a transitive $U(\boldd)$-space with stabilizer subgroup
$\aut_{\la}^u(U)$.  Since each fibre of the orbit map is obtainable as the
solution set of a system of linear equations, this map is separable, which
shows that, as a left
$U(\boldd)$-space,  $O(x)$ is isomorphic to the homogeneous space $U(\boldd) /
\aut_{\la}^u(U)$.  Note that the latter is the geometric quotient relative
to the action of $\aut_{\la}^u(U)$ on
$U(\boldd)$ by right translations.  On the other hand, the isomorphism given in
Lemma 3 induces an isomorphism $H \times (O(x) \cap V_p) \cong O(x)$ which
is equivariant with respect to  left multiplication by elements of $H$. Thus
$O(x) \cap V_p$ is the geometric quotient of
$O(x)$ by the left action of $H$.  Combining these facts, \cite{\Bong, Lemma
5.9(b)} now yields that $O(x) \cap V_p$ is the geometric quotient of $U(\boldd)$
modulo the action of the product $H\times
\aut_{\la}^u(U)$ by left resp. right translations.

To verify that this latter geometric quotient is isomorphic to $A / B$, we
will re-obtain  it along an alternate route.  Namely, we will again take two
successive quotients of $U(\boldd)$, but this time we will first divide by
$H$ and then by $\aut_{\la}^u(U)$.  Observe that the isomorphism $H \times A
\rightarrow U(\boldd)$ of Lemma 4 is $H \times \aut_{\la}^u(U)$-equivariant: we
keep the above action on
$U(\boldd)$, let
$H$ act on $H \times \aut_{\la}^u(U)$ from the left in the obvious fashion,
and define the right action of $\aut_{\la}^u(U)$ on $H \times
\aut_{\la}^u(U)$ as follows:  $(h,a).g = (h\cdot ga(g),a\cdot a(g)^{-1})$; 
in checking that this is well-defined, keep in mind that the subgroups $A$
and $\aut_{\la}^u(U)$ of $U(\boldd)$ commute, and note that, for $g \in
\aut_{\la}^u(U)$, we have $a(g) = (g_1,id)^{-1}$ in the notation of part
(1).  Consequently, the geometric quotient of
$U(\boldd)$ by the left action of $H$ is isomorphic to $A$, and the right
$\aut_{\la}^u(U)$-action on this partial quotient boils down to right
multiplication by elements of $B$.  We conclude that the left
$A$-space  $U(\boldd) / (H \times \aut_{\la}^u(U))$ is isomorphic to the
geometric quotient $A / B$ of
$A$ modulo the right action of $B$ by translations as required.  This
completes our argument.
\qed \enddemo 

In light of Theorem A$'$, the problem of recognizing the endomorphism algebras of
uniserial $\la$-modules among the finite dimensional commutative local
$K$-algebras imposes itself.  We leave this point to a subsequent investigation.

\enddocument